\documentclass[11pt]{article}
\usepackage{color}
\usepackage{tikz}
\usepackage{amssymb, amsmath}
\usepackage{epsfig}
\usepackage{graphicx}

\newtheorem{theorem}{Theorem}
 \newtheorem{corollary}{Corollary}
 \newtheorem{proposition}{Proposition}

\newcommand{\proof}{\noindent\textbf{Proof.~}}
\newcommand{\rr}{\mathbb{R}}
\newcommand{\vep}{\varepsilon}
\newcommand{\qed}{\space\hfill\hspace*{\fill} $\vbox{\hrule\hbox{\vrule
height1.3ex\hskip1.3ex\vrule}\hrule}$\hss\vskip\topsep\relax}

\begin{document}

\title{Egodic Theorems  for  discrete Markov chains }

\author{Nikolaos Halidias \\
{\small\textsl{Department of Mathematics }}\\
{\small\textsl{University of the Aegean }}\\
{\small\textsl{Karlovassi  83200  Samos, Greece} }\\
{\small\textsl{email: nikoshalidias@hotmail.com}}}

\maketitle

\begin{abstract} Let $X_n$ be a discrete time Markov chain with
state space $S$ (countably infinite, in general) and initial
probability distribution $\mu^{(0)} =
(P(X_0=i_1),P(X_0=i_2),\cdots,)$. What is the probability of
choosing in random some $k \in \mathbb{N}$ with $k \leq n$  such
that $X_k = j$ where $j \in S?$ This probability is the average
$\frac{1}{n} \sum_{k=1}^n \mu^{(k)}_j$ where $\mu^{(k)}_j = P(X_k
= j)$. In this note we will study the limit of this average
without assuming  that the chain is irreducible, using elementary
mathematical tools. Finally, we study the limit of the average
$\frac{1}{n} \sum_{k=1}^n g(X_k)$ where $g$ is a given function
for a  Markov chain not necessarily irreducible.
\end{abstract}



\section{Introduction}

Let $X_n$ be a discrete time Markov chain with state space $S$
(countably infinite, in general) and initial probability
distribution $\mu^{(0)}$, that is $\mu^{(0)}_i = P(X_0 = i)$ where
$i \in S$.  We will study the limit of the average
\begin{eqnarray*}
 \frac{1}{n} \sum_{k=1}^n \mu^{(k)}_j
\end{eqnarray*}
This quantity gives the probability of choosing in random an
integer $k$ with $ k \leq n$  such that $X_k = j$. Note that, for
any $i,j \in S$, we have
\begin{eqnarray}\label{nonumber}
\frac{1}{n} \sum_{k=1}^n \mu^{(k)}_j & = & \frac{1}{n}
\sum_{k=1}^n (\mu^{(0)} \cdot P^k)_j \nonumber \\  & = &
\frac{1}{n} \sum_{k=1}^n \sum_{i \in S} \mu^{(0)}_i P^k_{ij} \nonumber\\
 & = & \sum_{i \in S} \mu^{(0)}_i \frac{1}{n}
\sum_{k=1}^n P^k_{ij}
\end{eqnarray}
Therefore, one can study the desired limit by studying  the limit
of the average $\frac{1}{n} \sum_{k=1}^n P^k_{ij}$. To do so  one
can  use the limit theorems for $P^n_{ij}$ (see for example
\cite{Ash}) and the well known fact that if $a_n \to a$ then
$\frac{1}{n} \sum_{k=1}^n a_k \to a$. However, here we will give a
different proof without using the limit theorems and without
assuming that the chain is irreducible. Moreover, we will study
the behavior of the limit $\frac{1}{n} \sum_{k=1}^n g(X_k)$ for a
given function $g$, using elementary mathematical tools.

\section{The main results}
Let $X_n$ be a Markov chain with (countably infinite in general)
state space $S$.

\begin{theorem} It holds that, for any $i,j \in S$,
\begin{eqnarray*}
\lim_{n \to \infty} \frac{1}{n} \sum_{k=1}^{n} \mu_j^{(k)} = \left\{%
\begin{array}{cl}
   \displaystyle  \frac{1}{m_j} \sum_{i \in S} \mu^{(0)}_i f_{ij}, & \mbox{ when } j \mbox{ is positive recurrent }
    \\[0.5cm]
    0, & \hbox{ otherwise } \\
\end{array}%
\right.
\end{eqnarray*}
and
\begin{eqnarray*}
\lim_{n \to \infty} \frac{ \sum_{k=1}^{n} P^k_{ij}}{n} = \left\{%
\begin{array}{ll}
    \frac{f_{ij}}{m_j}, & \mbox{ when } j \mbox{ is positive recurrent  } \\[0.2cm]
    0, & \hbox{ otherwise } \\
\end{array}%
\right.
\end{eqnarray*}
 where $f_{ij} = P(\exists \; n \in \mathbb{N} : X_n = j|X_0
=i)$.
\end{theorem}
\proof We know (see \cite{Ash}) that when $j$ is transient or null
recurrent $\lim_{n \to \infty} P^n_{ij} = 0$. Therefore $\lim_{n
\to \infty} \frac{1}{n} \sum_{k=1}^n P^n_{ij} =0$ and using
\ref{nonumber} the result follows. Next we suppose that $j$ is
positive recurrent.

Let the random variables $N_j^k = \left\{%
\begin{array}{ll}
    1, & \mbox{ when } X_k = j \\
    0, & \mbox{ otherwise} \\
\end{array}%
\right.$ and $M_j(n) = \sum_{k=1}^n N_j^k$. Because
\begin{eqnarray}\label{EM_jn_n}
\mathbb{E} \left( \frac{M_j(n)}{n} \right) = \frac{1}{n}
\sum_{k=1}^n \mathbb{E} N_j^k = \frac{1}{n} \sum_{k=1}^n P(X_k =
j) =\frac{1}{n} \sum_{k=1}^n \mu^{(k)}_j
\end{eqnarray}
we will study the quantity $\mathbb{E} \left( \frac{M_j(n)}{n}
\right)$.

Let the  event   $A_i = \{ \exists \; n \in \mathbb{N} : X_n = j
\} \cap \{ X_0 = i \} $ where $i \in S$. Because $P(A_i) =
P(\exists \; n \in \mathbb{N} : X_n = j|X_0 =i) \cdot \mu^{(0)}_i$
we see that  $P(A_i) = f_{ij} \cdot \mu^{(0)}_i$ where $f_{ij} =
P(\exists \; n \in \mathbb{N} : X_n = j|X_0 =i)$.

We will work under the probability measure
 $P_{A_i}(\cdot) = P( \cdot |A_i)$ while the corresponding expected value will be denoted by  $\mathbb{E}_{A_i}$.

We define the following sequence of random variables,

\begin{eqnarray*}
 n_1(\omega) & = & \left\{
\begin{array}{ll}
    \min \{ n \in \mathbb{N} : X_n(\omega) = j \}, & \mbox{ when } \omega \in A_i \\[0.1cm]
    \infty, & \mbox{ otherwise } \\
\end{array}%
\right. \\[0.3cm]
n_2(\omega) & = & \left\{
\begin{array}{ll}
    \min \{ n > n_1 : X_n(\omega) = j \}, & \mbox{ when } \omega \in A_i \\[0.1cm]
    \infty, & \mbox{ otherwise } \\
\end{array}%
\right. \\[0.3cm]
& \vdots & \\[0.3cm]
n_k(\omega) & = & \left\{
\begin{array}{ll}
    \min \{ n > n_{k-1}: X_n(\omega) = j \}, & \mbox{ when } \omega \in A_i \\[0.1cm]
    \infty, & \mbox{ otherwise } \\
\end{array}%
\right.
\end{eqnarray*}
  We define also  $Z_m = \left\{%
\begin{array}{ll}
    n_{m+1} - n_m, & \mbox{ when } \omega \in A_i \\
    0, & \mbox{ otherwise } \\
\end{array}%
\right.    $ for $m \geq 1$ which gives us the number of
transitions needed to return back to  $j$. Note that the sequence
 $Z_1, Z_2,\cdots,$ is an independent and identically distributed  sequence of
 random variables. The mean recurrent time $m_j$ is such that
 $m_j = \mathbb{E}_{A_i}(Z_k)$
for every $k \geq 1$. Next we define the random variable
 $S_l = Z_1 + \cdots +
Z_l$ with $S_0 = 0$. Note that
\begin{eqnarray}\label{s_ln_1}
S_l + n_1 = n_{l+1} \quad \mbox{ for every } l \geq 0
\end{eqnarray}
Using the strong law of large numbers we have that
\begin{eqnarray*}
P_{A_i} \left(\{\omega \in \Omega: \lim_{n \to \infty}
\frac{S_n}{n} = m_j \}\right) = 1
\end{eqnarray*}

Note that  $M_j(n) \to \infty$ as $n \to \infty$ for almost all
$\omega \in \Omega$ when $j$ is
 recurrent and its easy to see that $n_{M_j(k)} \leq k$
for every $k \geq 1$.

Using \ref{s_ln_1} we see that the following inequality hold
\begin{eqnarray*}
S_{M_j(n)-1}+n_1 \leq n \leq S_{M_j(n)}+n_1, \quad n \geq 1, \quad
\mbox{ for every } \omega \in A_i
\end{eqnarray*}
Dividing by $M_j(n) > 0$ for $n > n_1$ we get
\begin{eqnarray*}
\frac{S_{M_j(n)-1}+n_1}{M_j(n)-1} \frac{M_j(n)-1}{M_j(n)} \leq
\frac{n}{M_j(n)}  \leq \frac{S_{M_j(n)}}{M_j(n)}, \quad n \geq
n_2, \quad \mbox{ for every } \omega \in A_i
\end{eqnarray*}
 Therefore it holds that
\begin{eqnarray}\label{pithanotitasigklisis}
P_{A_i} \left( \{\omega \in \Omega: \lim_{n \to \infty}
\frac{M_j(n)}{n} = \frac{1}{m_j} \} \right) = 1
\end{eqnarray}

Next we will study the limit of the quantity
\begin{eqnarray*}
\lim_{n \to \infty} \frac{\mathbb{E}_{A_i}(M_j(n) )}{n}
\end{eqnarray*}
Using the dominated convergence theorem it follows that
\begin{eqnarray*}
\lim_{n \to \infty} \frac{\mathbb{E}_{A_i}(M_j(n) )}{n} & = &
\mathbb{E}_{A_i}\left(\lim_{n \to \infty} \frac{M_j(n)}{n} \right)
\\ & = &  \mathbb{E}_{A_i} \left( \frac{1}{m_j}  \right)  \\ & = &
\frac{1}{m_j}
\end{eqnarray*}
But, since
\begin{eqnarray*}
\mathbb{E}_{A_i} \left( \frac{M_j(n)}{n}   \right) =
\frac{\mathbb{E} \left( \frac{M_j(n)}{n} \mathbb{I}_{A_i} \right)
}{P(A_i)}
\end{eqnarray*}
 it follows that
\begin{eqnarray}\label{m0fijmj}
\lim_{n \to \infty} \mathbb{E} \left( \frac{M_j(n)}{n}
\mathbb{I}_{A_i}  \right) = \frac{P(A_i)}{m_j} = \mu^{(0)}_i
\frac{f_{ij}}{m_j}
\end{eqnarray}

 Because
\begin{eqnarray*}
\mathbb{E} \left(\frac{M_j(n)}{n}\right) = \sum_{ i \in S}
\mathbb{E} \left(\frac{M_j(n)}{n} \mathbb{I}_{A_i} \right)
\end{eqnarray*}
 we  obtain using  \ref{m0fijmj}
\begin{eqnarray*}
\lim_{n \to \infty} \mathbb{E} \left(\frac{M_j(n)}{n}\right) & = &
\lim_{n \to \infty}\sum_{ i \in S} \mathbb{E}
\left(\frac{M_j(n)}{n} \mathbb{I}_{A_i} \right) \\ & = & \sum_{i
\in S} \lim_{n \to \infty} \mathbb{E} \left(\frac{M_j(n)}{n}
\mathbb{I}_{A_i} \right) \\ & = & \sum_{i \in S} \mu^{(0)}_i
\frac{f_{ij}}{m_j} \\ & = & \frac{1}{m_j} \sum_{i \in S}
\mu^{(0)}_i f_{ij}
\end{eqnarray*}
where we have used the dominated convergence theorem to get the
second equality above.

If $m_{ij}(n) = \mathbb{E}(M_j(n)|X_0=i)$ then we have
\begin{eqnarray*}
m_{ij}(n) & = & \mathbb{E} \left(M_j(n) | X_0  =  i \right)
\\ & = & \mathbb{E} \left( \sum_{k=1}^{n} N_j^k | X_0 = i \right) \\ &
=& \sum_{k=1}^{n} \mathbb{E} (N_j^k | X_0 = i) \\ & = &
\sum_{k=1}^{n} P_{ij}^k
\end{eqnarray*}

 Denoting by $A = \{ \exists \; k \in \mathbb{N} : X_k = j \}$,
we have
\begin{eqnarray*}
\mathbb{E}\left(\frac{M_j(n)}{n} |X_0 = i\right) & = &
\mathbb{E}\left(\frac{M_j(n)}{n} \mathbb{I}_A |X_0 = i\right) +
\mathbb{E} \left(\frac{M_j(n)}{n} \mathbb{I}_{A^c} |X_0 = i\right)
\\ &  = & \mathbb{E}\left(\frac{M_j(n)}{n} \mathbb{I}_A |X_0 =
i\right)\\ & = & \frac{ \mathbb{E} \left( \frac{M_j(n)}{n}
\mathbb{I}_{A_i} \right)} { m^{(0)}_i}
\end{eqnarray*}
because $M_j(n) \mathbb{I}_{A^c} = 0$. That means that  $$\lim_{n
\to \infty}\frac{m_{ij}(n)}{n} = \lim_{n \to \infty}
\mathbb{E}\left(\frac{M_j(n)}{n} |X_0 = i\right) =
 \frac{f_{ij}}{m_j}$$
Therefore
\begin{eqnarray*}
\lim_{n \to \infty} \frac{ \sum_{k=1}^{n} P^k_{ij}}{n} = \left\{%
\begin{array}{ll}
    \frac{f_{ij}}{m_j}, & \mbox{ when } j \mbox{ is positive recurrent  } \\[0.2cm]
    0, & \hbox{ otherwise } \\
\end{array}%
\right.
\end{eqnarray*}
 \qed

\begin{proposition}\label{almost surely} It holds that, when $j$ is positive recurrent,
\begin{eqnarray*}
\{ \omega \in \Omega: \lim_{n \to \infty} \frac{M_j(n)}{n} =
\frac{1}{m_j} \} \cup \{ \omega \in \Omega: \lim_{n \to \infty}
\frac{M_j(n)}{n} = 0 \} = \Omega \smallsetminus E
\end{eqnarray*}
with $P(E) = 0$. More precisely, it holds that
\begin{eqnarray*}
P\left( \{ \omega \in \Omega: \lim_{n \to \infty} \frac{M_j(n)}{n}
= \frac{1}{m_j}\} \right)=\sum_{i \in S} \mu_i^{(0)} \cdot f_{ij}
\end{eqnarray*}
and
\begin{eqnarray*}
P\left( \{ \omega \in \Omega: \lim_{n \to \infty} \frac{M_j(n)}{n}
= 0 \} \right)=\sum_{i \in S} \mu_i^{(0)} \cdot (1-f_{ij})
\end{eqnarray*}
If $j$ is null recurrent or transient, then
\begin{eqnarray*}
P\left( \{ \omega \in \Omega: \lim_{n \to \infty} \frac{M_j(n)}{n}
= 0 \} \right)=1
\end{eqnarray*}
\end{proposition}
\proof

$\cdot$ Assume that $j$ is positive recurrent. Denoting by $B = \{
\omega \in \Omega: \lim_{n \to \infty} \frac{M_j(n)}{n} =
\frac{1}{m_j} \}$ we can write
\begin{eqnarray*}
B = \bigcup_{i \in S} \{ \omega \in \Omega: \lim_{n \to \infty}
\frac{M_j(n)}{n} = \frac{1}{m_j} \} \cap \{ X_0 = i \} =
\bigcup_{i \in S} B_i
\end{eqnarray*}
and therefore $P(B) = \sum_{i \in S} P(B_i)$.

But
\begin{eqnarray*}
B_i = B_i \cap \{ \exists \; k\in \mathbb{N} : X_k = j \}  \bigcup
B_i \cap \{ \nexists \; k\in \mathbb{N} : X_k = j \}
\end{eqnarray*}
so $P(B_i) = P(B_i \cap \{ \exists \; k\in \mathbb{N} : X_k = j \}
 ) + P(B_i \cap \{ \nexists \; k\in \mathbb{N} : X_k = j \}
)$. Recalling \ref{pithanotitasigklisis} we can write that
\begin{eqnarray*}
P(B_i \cap \{ \exists \; k\in \mathbb{N} : X_k = j \}
 ) = P_{A_i} \left(\{ \omega \in \Omega: \lim_{n \to \infty} \frac{M_j(n)}{n}
= \frac{1}{m_j}\} \right) \cdot P(A_i) = \mu_i^{(0)} \cdot f_{ij}
\end{eqnarray*}
Moreover
\begin{eqnarray*}
P(B_i \cap \{ \nexists \; k\in \mathbb{N} : X_k = j \}) = 0
\end{eqnarray*}
since in this event $M_j(n) = 0$. Therefore $P(B_i) = \mu_i^{(0)}
\cdot f_{ij}$ and thus $$P(B) = \sum_{i \in S} \mu_i^{(0)} \cdot
f_{ij}$$

Denote now $\Gamma_i = \{ \nexists \; k\in \mathbb{N} : X_k = j
\}\cap \{X_0= i \} $. Then
\begin{eqnarray*}
P_{\Gamma_i} \left( \{ \omega \in \Omega: \lim_{n \to \infty}
\frac{M_j(n)}{n} = 0 \} \right) = 1
\end{eqnarray*}
where $P_{\Gamma_i}(\cdot) = P(\cdot | \Gamma_i)$. Thus
\begin{eqnarray*}
P \left( \{ \omega \in \Omega: \lim_{n \to \infty}
\frac{M_j(n)}{n} = 0 \} \cap \Gamma_i \right) = P(\Gamma_i) =
\mu_i^{(0)} (1-f_{ij})
\end{eqnarray*}
That means that
\begin{eqnarray*}
P \left( \{ \omega \in \Omega: \lim_{n \to \infty}
\frac{M_j(n)}{n} = 0 \} \cap \Gamma \right) = \sum_{i \in S}
\mu_i^{(0)} (1-f_{ij})
\end{eqnarray*}
where $\Gamma = \{ \nexists \; k\in \mathbb{N} : X_k = j \}$. Thus
\begin{eqnarray*}
P \left( \{ \omega \in \Omega: \lim_{n \to \infty}
\frac{M_j(n)}{n} = 0 \}  \right) \geq \sum_{i \in S} \mu_i^{(0)}
(1-f_{ij})
\end{eqnarray*}
The events
\begin{eqnarray*}
\{ \omega \in \Omega: \lim_{n \to \infty} \frac{M_j(n)}{n} =
\frac{1}{m_j} \} \quad \mbox{ και } \quad \{ \omega \in \Omega:
\lim_{n \to \infty} \frac{M_j(n)}{n} = 0 \}
\end{eqnarray*}
are disjoint, therefore
\begin{eqnarray*}
1 & \leq & P \left(\{ \omega \in \Omega: \lim_{n \to \infty}
\frac{M_j(n)}{n} = \frac{1}{m_j} \} \right) + P\left(\{ \omega \in
\Omega: \lim_{n \to \infty} \frac{M_j(n)}{n} = 0 \} \right) \\ & =
&  P \left(\{ \omega \in \Omega: \lim_{n \to \infty}
\frac{M_j(n)}{n} = \frac{1}{m_j} \} \cup \{ \omega \in \Omega:
\lim_{n \to \infty} \frac{M_j(n)}{n} = 0 \} \right)\\ & \leq   & 1
\end{eqnarray*}
Therefore
\begin{eqnarray*}
\{ \omega \in \Omega: \lim_{n \to \infty} \frac{M_j(n)}{n} =
\frac{1}{m_j} \} \cup \{ \omega \in \Omega: \lim_{n \to \infty}
\frac{M_j(n)}{n} = 0 \} = \Omega \smallsetminus E
\end{eqnarray*}
with $P(E) = 0$ and
\begin{eqnarray*}
P \left( \{ \omega \in \Omega: \lim_{n \to \infty}
\frac{M_j(n)}{n} = 0 \}  \right) = \sum_{i \in S} \mu_i^{(0)}
(1-f_{ij})
\end{eqnarray*}

$\cdot$ Assume now that $j$ is null recurrent
 and let the sequence of random variables
 $Z_m = \left\{%
\begin{array}{ll}
    n_{m+1} - n_m, & \mbox{ όταν } \omega \in A_i \\
    0, & \mbox{ otherwise } \\
\end{array}%
\right.    $ for $m \geq 1$. Because $j$ is null recurrent we have
that  $\mathbb{E}(Z_m) = \infty$ for every  $m \geq 1$. We define
now the sequence  $Z_m^R = Z_m \mathbb{I}_{ \{ Z_m < R \} }$ for
$R
> 0$ for which it holds that  $\mathbb{E}(Z_m^R) < \infty$ for every  $m \geq 1$. Moreover,  $\mathbb{E}(Z_1^R) = \mathbb{E}(Z_m^R)$
for every $m \geq 1$.  This sequence   is again an independent and
identical distributed sequence of random variables. Therefore we
can use the strong law of large numbers to get
\begin{eqnarray*}
P_{A_i} \left( \omega \in \Omega : \lim_{n \to \infty}
\frac{S_n^R}{n} = \mathbb{E}_{A_i}(Z_1^R) \right) = 1
\end{eqnarray*}
where $S_n^R = Z_1^R + Z_2^R + \cdots + Z_n^R \leq S_n = Z_1 +
\cdots + Z_n$ and  $A_i$, $P_{A_i}$ is as before. Therefore it
holds that
\begin{eqnarray*}
S_{M_j(n)-1}^R + n_1 \leq S_{M_j(n)-1} + n_1 \leq n
\end{eqnarray*}
So
\begin{eqnarray*}
\frac{S^R_{M_j(n)-1}+n_1}{M_j(n)-1} \frac{M_j(n)-1}{M_j(n)} \leq
\frac{n}{M_j(n)}  , \quad n \geq n_2, \quad \mbox{ for every }
\omega \in A_i
\end{eqnarray*}
Letting   $n \to \infty$ we get that
\begin{eqnarray*}
0 \leq \limsup_{n \to \infty} \frac{M_j(n)}{n} \leq
\frac{1}{\mathbb{E}(Z_m^R)}, \quad \mbox{ almost surely, } \;
\mbox{ for every  } R> 0
\end{eqnarray*}
under the probability measure  $P_{A_i}$.  Note that  $Z_m^R$ is
an increasing sequence in  $R$ and that  $Z_m^R \to Z_m$ as $R \to
\infty$ almost surely. Therefore  $\mathbb{E}_{A_i} (Z_m^R) \to
\mathbb{E}_{A_i} (Z_m) = \infty$ using the monotone convergence
theorem. That means that
\begin{eqnarray*}
\lim_{n \to \infty} \frac{M_j(n)}{n} = 0  \quad \mbox{ almost
surely }
\end{eqnarray*}
under the probability measure $P_{A_i}$, i.e.
\begin{eqnarray}\label{p_a_j=1}
P_{A_i} \left( \{ \omega \in \Omega : \lim_{n \to \infty}
\frac{M_j(n)}{n} =0 \} \right) =1
\end{eqnarray}

Let now  the event $\{ \omega \in \Omega : \limsup_{n \to \infty}
\frac{M_j(n)}{n} \geq \vep \}$ where $\vep > 0$. Noting that
\begin{eqnarray*}
 P \left( \{ \omega
\in \Omega : \limsup_{n \to \infty} \frac{M_j(n)}{n} \geq \vep \}
\cap A^c \right) = 0
\end{eqnarray*}
where $A = \{ \exists \; l \in \mathbb{N}  :X_l = j \}$ and
\begin{eqnarray*}
P \left( \{ \omega \in \Omega : \limsup_{n \to \infty}
\frac{M_j(n)}{n} \geq \vep \} \cap A \right) & = &  \sum_{i \in S}
P \left(\{ \omega \in \Omega : \limsup_{n \to \infty}
\frac{M_j(n)}{n} \geq \vep \} \cap A_i \right) \\ & = & \sum_{i
\in S} \underbrace{P_{A_i} \left( \{ \omega \in \Omega :
\limsup_{n \to \infty} \frac{M_j(n)}{n} \geq \vep \} \right)}_{ =
0, \mbox{ see \ref{p_a_j=1}}}  P(A_i) \\ & = & 0
\end{eqnarray*}
we obtain
\begin{eqnarray*}
P\left( \{ \omega \in \Omega : \limsup_{n \to \infty}
\frac{M_j(n)}{n} \geq \vep \} \right) = 0
\end{eqnarray*}
Because $\frac{M_j(n)}{n} \geq 0$  it follows the desired result.

 $\cdot$ Finally we assume that $j$ is transient.
It is well known that ι $P(M_j < \infty |X_0 = i) = 1$ for every
state  $i$, where $M_j = \lim_{n \to \infty} M_j(n)$. Therefore
\begin{eqnarray*}
P(M_j < \infty) = \sum_{i \in S} P(M_j < \infty |X_0 = i) \cdot
P(X_0=i) = \sum_{i \in S} \mu^{(0)}_i = 1
\end{eqnarray*}

Moreover
\begin{eqnarray*}
\Omega = \left( \bigcup_{N=0}^{\infty} B_N \right) \cup B_{\infty}
\end{eqnarray*}
where $B_N = \{ M_j = N \}$ and $B_{\infty} = \{ M_j = \infty \}$.
Thus
\begin{eqnarray*}
\sum_{N=0}^{\infty} P(B_N) = 1
\end{eqnarray*}
since $P(B_{\infty}) = 0$.

Therefore we can write
\begin{eqnarray*}
& & \{ \omega \in \Omega : \lim_{n \to \infty} \frac{M_j(n)}{n} =
0 \} \\ & = &  \left( \bigcup_{N=0}^{\infty} \{ \omega \in \Omega
: \lim_{n \to \infty} \frac{M_j(n)}{n} = 0 \} \cap B_N \right)
\cup \{ \omega \in \Omega : \lim_{n \to \infty} \frac{M_j(n)}{n} =
0 \} \cap B_{\infty}
\end{eqnarray*}

Thus
\begin{eqnarray*}
P \left( \{ \omega \in \Omega : \lim_{n \to \infty}
\frac{M_j(n)}{n} = 0 \}  \right) = \sum_{N=0}^{\infty} P \left(
 \{ \omega \in \Omega : \lim_{n \to \infty}
\frac{M_j(n)}{n} = 0 \} \cap B_N \right)
\end{eqnarray*}
since $P\left(\{ \omega \in \Omega : \lim_{n \to \infty}
\frac{M_j(n)}{n} = 0 \} \cap B_{\infty} \right) \leq P(B_{\infty})
= 0$. But
\begin{eqnarray*}
& & P \left(
 \{ \omega \in \Omega : \lim_{n \to \infty}
\frac{M_j(n)}{n} = 0 \} \cap B_N \right) \\ & = & P \left( \{
\omega \in \Omega : \lim_{n \to \infty} \frac{M_j(n)}{n} = 0 \} |
B_N \right) P(B_N) \\ & = & P(B_N)
\end{eqnarray*}
since it holds that  $P \left( \{ \omega \in \Omega : \lim_{n \to
\infty} \frac{M_j(n)}{n} = 0 \}   |  B_N \right) =1$. Since
$\sum_{N=0}^{\infty} P(B_N) = 1$ we obtain the desired result.
 \qed

\begin{corollary}
If $g:S \to \rr$ is such that
\begin{eqnarray*}
\sum_{i \in S}  |g(i)| < \infty
\end{eqnarray*}
then it holds that
\begin{eqnarray*}
\lim_{n \to \infty} \frac{1}{n} \sum_{k=1}^{n} \mathbb{E} g(X_k) =
\sum_{j \in C} \frac{g(j)}{m_j} \sum_{i \in S} \mu^{(0)}_i f_{ij}
\end{eqnarray*}
where $C \subseteq S$ is the subset of $S$ of positive recurrent
states.
\end{corollary}
\proof Note that $\displaystyle g(X_k) = \sum_{j \in S} g(j)
\mathbb{I}_{ \{ X_k = j \} }$. Therefore
\begin{eqnarray*}
\frac{1}{n} \sum_{k=1}^n \mathbb{E} g(X_k) & = & \frac{1}{n}
\sum_{k=1}^n \sum_{j \in S} g(j) \mathbb{E} \mathbb{I}_{ \{ X_k =
j \} } \\ & = & \sum_{j \in S} g(j) \frac{1}{n} \sum_{k=1}^n
\mathbb{E} \mathbb{I}_{ \{ X_k = j \} }
\\ & = & \sum_{j \in S} g(j)  \mathbb{E} \left(\frac{M_j(n)}{n}
\right)
\end{eqnarray*}
We have interchange the sums $\sum_{k=1}^n \sum_{j \in S}$ because
the series is absolutely convergent since $\sum_{i \in S}  |g(i)|
< \infty$.

 So
\begin{eqnarray*}
\lim_{n \to \infty} \frac{1}{n} \sum_{k=1}^n \mathbb{E}g(X_k) =
\lim_{n \to \infty} \sum_{j \in S} g(j) \mathbb{E}
\left(\frac{M_j(n)}{n} \right) = \sum_{j \in C} \frac{g(j)}{m_j}
\sum_{i \in S} \mu^{(0)}_i f_{ij}
\end{eqnarray*}
 We have used the dominated convergence theorem to interchange
the limit with the sum in the second equality above.
 \qed

\begin{corollary}\label{corobirkoff} Given a function $g:S \to \rr$ such that
\begin{eqnarray*}
\sum_{i \in S}  |g(i)| < \infty
\end{eqnarray*}
it holds that
\begin{eqnarray*}
\lim_{n \to \infty} \frac{1}{n} \sum_{k=1}^n g(X_k) = \sum_{j \in
C} \frac{g(j)}{m_j} \mathbb{I}_{A^j} \quad \mbox{ almost surely}
\end{eqnarray*}
where $A^j = \{\omega \in \Omega: \exists \; l \in \mathbb{N} :
X_l = j \}$.
\end{corollary}
\proof Note that
\begin{eqnarray*}
\frac{1}{n} \sum_{k=1}^n  g(X_k) & = & \frac{1}{n} \sum_{k=1}^n
\sum_{j \in S} g(j) \mathbb{I}_{ \{ X_k = j \}} \\ & = & \sum_{j
\in S} g(j) \frac{M_j(n)}{n} \\ & = & \sum_{j \in C} g(j)
\frac{M_j(n)}{n} \mathbb{I}_{A^j} + \sum_{j \in NR} g(j)
\frac{M_j(n)}{n} + \sum_{j \in T} g(j) \frac{M_j(n)}{n}
\end{eqnarray*}
where $C \subseteq S$ is the subset of positive recurrent states
of $S$, $NR \subseteq S$ is the subset of null recurrent states of
$S$, $T \subseteq S$ is the subset of transient states of $S$ and
$A^j = \{\omega \in \Omega: \exists \; l \in \mathbb{N} : X_l = j
\}$. The condition on $g$, i.e. $\sum_{i \in S} |g(i)| < \infty$
is needed in order to interchange the sums to get the second
equation above.

Note that $A^j = \bigcup_{i \in S} \{ \exists \; l \in \mathbb{N}
:X_l = j \} \cap \{X_0=i \} $ and therefore $P(A^j) = \sum_{i \in
S} \mu_i^{(0)} \cdot f_{ij}$.

 Finally, using
proposition \ref{almost surely}, we obtain the desired result,
i.e.
\begin{eqnarray*}
\lim_{n \to \infty} \frac{1}{n} \sum_{k=1}^n  g(X_k) = \sum_{j \in
C} \frac{g(j)}{m_j}\mathbb{I}_{A^j}, \quad \mbox{ almost surely }
\end{eqnarray*}
 \qed
This corollary is closely related to Birkoff's ergodic theorem
(see for example \cite{Durrett} and  \cite{Skorohod}).


\begin{thebibliography}{99}


\bibitem{Ash1} R. Ash, {\em Basic Probability Theory}, Wiley,
1970.

 \bibitem{Ash} R. Ash - C.
Doleans-Dade, {\em Probability and Measure Theory}, Elsevier,
2000.


\bibitem{Bremaunt} P. Bremaud, {\em Markov Chains, Gibbs Fields, Monte Carlo Simulation and Queues}, Springer, 1999.

\bibitem{Durrett} R. Durrett, {\em Probabiblity: Theory and
Examples}, Cambridge University Press, 2010.


\bibitem{Grimmett} G. Grimmett - D. Stirzaker, {\em Probability
and Random Processes}, Oxford University Press, 2001.

\bibitem{Hoel} P. Hoel -  S. Port - C. Stone, {\em Introduction to
Stochastic Processes}, Houghton Mifflin Company, 1972.


\bibitem{Skorohod} I. Gikhman - A. Skorokhod, {\em Introduction
to the Theory of Random Processes}, 1969.

\bibitem{Jacod} J. Jacod - P. Protter, {\em Probability Essentials}, Springer,
2004.

\bibitem{Karlin} S. Karlin - H. Taylor, {\em A First Course in
Stochastic Processes}, Academic Press, 1975.

\bibitem{Kemeny} J. Kemeny - J. Snell, {\em Finite Markov Chains},
Springer, 1976.

\bibitem{Norris} J. Norris, {\em Markov Chains}, Cambridge
University Press, 1997.

\bibitem{Revuz} D. Revuz, { \em Markov Chains}, North Holland,
1984.


\bibitem{Stroock} D. Stroock, {\em An Introduction to Markov
Processes}, Springer, 1999.


\end{thebibliography}
\end{document}